\documentclass[12pt]{scrartcl}

\usepackage[
  pdftex,
  bookmarks,
  bookmarksopen=true,
  bookmarksnumbered=true,
  pdfusetitle,
  colorlinks,
  plainpages=false,
  ]{hyperref}

\usepackage{amsfonts}
\usepackage{amsmath}
\usepackage{amssymb}
\usepackage{amsthm}
\usepackage{graphicx}
\usepackage{color}
\usepackage{sidecap}
\usepackage{wasysym}
\usepackage{prettyref}
\usepackage{mathrsfs}
\usepackage[utf8]{inputenc}

\newcommand{\N}{\ensuremath{\mathbb{N}_0}}

\renewcommand{\S}{\ensuremath{\mathbb{S}}}

\newcommand{\R}{\ensuremath{\mathbb{R}}}
\newcommand{\C}{\ensuremath{\mathbb{C}}}

\newcommand{\abs}[1]{\ensuremath{\left\vert#1\right\vert}}
\newcommand{\dx}{\mathrm{d}}

\newcommand{\zb}[1]{\ensuremath{\boldsymbol{#1}}}

\renewcommand{\d}{\, \mathrm{d}}

\newcommand{\norm}[1]{\left\lVert #1
	\right\rVert}

\newtheorem{theorem}{Theorem}[section]
\newtheorem{lemma}[theorem]{Lemma}
\newtheorem{remark}[theorem]{Remark}
\newtheorem{definition}[theorem]{Definition}
\newtheorem{example}[theorem]{Example}
\newtheorem{corollary}[theorem]{Corollary}
\newtheorem{proposition}[theorem]{Proposition}
\newtheorem{problem}[theorem]{Problem}

\newtheorem*{theorem*}{Theorem}
\newtheorem*{lemma*}{Lemma}
\newtheorem*{definition*}{Definition}
\newtheorem*{example*}{Example}
\newtheorem*{corollary*}{Corollary}
\newtheorem*{proposition*}{Proposition}

\theoremstyle{remark}
\newtheorem*{remark*}{Remark}

\newenvironment{Theorem}[1][noisnotdefined]{ \ifthenelse{\equal{#1}{noisnotdefined}}{\begin{theorem}}{\begin{theorem}[#1]}\normalfont\slshape}{\end{theorem}}
\newenvironment{Lemma}{ \begin{lemma}\normalfont\slshape}{\end{lemma}}
\newenvironment{Remark}[1][noisnotdefined]{ \ifthenelse{\equal{#1}{noisnotdefined}}{\begin{remark}}{\begin{remark}[#1]}\normalfont\rmfamily}{\bend\end{remark}}

\newenvironment{Corollary}{ \begin{corollary}\normalfont\rmfamily}{\end{corollary}}
\newenvironment{Proposition}{ \begin{proposition}\normalfont\rmfamily}{\end{proposition}}

\newenvironment{Theorem*}[1][noisnotdefined]{ \ifthenelse{\equal{#1}{noisnotdefined}}{\begin{theorem*}}{\begin{theorem*}[#1]}\normalfont\slshape}{\end{theorem*}}
\newenvironment{Lemma*}{ \begin{lemma*}\normalfont\slshape}{\end{lemma*}}
\newenvironment{Remark*}{ \begin{remark*}\normalfont\rmfamily}{\bend\end{remark*}}
\newenvironment{Example*}{ \begin{example*}\normalfont\rmfamily}{\bend\end{example*}}
\newenvironment{Definition*}{ \begin{definition*}\normalfont\rmfamily}{\end{definition*}}
\newenvironment{Corollary*}{ \begin{corollary*}\normalfont\rmfamily}{\end{corollary*}}
\newenvironment{Proposition*}{ \begin{proposition*}\normalfont\rmfamily}{\end{proposition*}}

\numberwithin{equation}{section}

\newcommand{\bend}{\hspace*{0ex} \hfill \hbox{\vrule height
    1.5ex\vbox{\hrule width 1.4ex \vskip 1.4ex\hrule  width 1.4ex}\vrule
    height 1.5ex}}

\long\def\symbolfootnote[#1]#2{\begingroup%
\def\thefootnote{\fnsymbol{footnote}}\footnote[#1]{#2}\endgroup}

\newrefformat{def}{Definition \ref{#1}}
\newrefformat{sub}{Subsection \ref{#1}}
\newrefformat{prop}{Proposition \ref{#1}}
\newrefformat{fig}{Figure \ref{#1}}
\newrefformat{rem}{Remark \ref{#1}}

\title{The Funk--Radon transform for hyperplane sections through a common point}

\date{}
\author{%
	Michael Quellmalz\thanks{Chemnitz University of Technology, Faculty of
		Mathematics, D-09107 Chemnitz, Germany. 
		\newline E-mail: \href{mailto:michael.quellmalz@mathematik.tu-chemnitz.de}{michael.quellmalz@mathematik.tu-chemnitz.de}
		}
}

\begin{document}

\maketitle

\begin{abstract}
The Funk--Radon transform, also known as the spherical Radon transform, assigns to a function on the sphere its mean values along all great circles. Since its invention by Paul Funk in 1911, the Funk--Radon transform has been generalized to other families of circles as well as to higher dimensions.

We are particularly interested in the following generalization: we consider the intersections of the sphere with hyperplanes containing a common point inside the sphere. If this point is the origin, this is the same as the aforementioned Funk--Radon transform. We give an injectivity result and a range characterization of this generalized Radon transform by finding a relation with the classical Funk--Radon transform.

\medskip

\textit{Math Subject Classifications.}
45Q05, 
44A12.  


\bigskip

\end{abstract}


\section{Introduction}

The reconstruction of a function from its integrals along certain submanifolds is a key task in various mathematical fields, including the modeling of imaging modalities, cf.\ \cite{NaWue00}.
One of the first works was in 1911 by Funk \cite{Fun11}.
He considered what later became known by the names Funk--Radon transform, Minkowski--Funk transform or spherical Radon transform.
Here, the function is defined on the unit sphere and we know its mean values along all great circles.
Another famous example is the Radon transform \cite{radon17}, where a function on the plane is assigned to its mean values along all lines.
Over the last more than 100 years, the Funk--Radon transform has been generalized to other families of circles as well as to higher dimensions.

In this article, we draw our attention to integrals along certain subspheres of the $(d-1)$-dimensional unit sphere $\S^{d-1} = \{\zb\xi\in\R^{d}\mid \norm{\zb\xi}=1 \}$.
Any subsphere of $\S^{d-1}$ is the intersection of the sphere with a hyperplane.
In particular, we consider the subspheres of $\S^{d-1}$ whose hyperplanes have the common point 
$
(0,\dots,0,z)^\top
$
strictly inside the sphere for $z\in[0,1)$.
We define the spherical transform $\mathcal U_z f$ of a function $f\in C(\S^{d-1})$ by
\begin{equation*}
\mathcal U_z f (\zb{\xi}) 
= \frac{1}{V(\mathscr C_z^{\zb\xi})} 
\int_{\mathscr C_z^{\zb\xi}} f \d \mathscr C_z^{\zb\xi}
,\qquad \zb{\xi} \in\S^{d-1},
\end{equation*}
which computes the mean values of $f$ along the subspheres 
\begin{equation*}
\mathscr C_z^{\zb\xi}
= \{\zb\eta\in\S^{d-1}\mid \langle\zb\eta,\zb\xi\rangle=z\xi_{d} \}
,\qquad z\in [0,1),
\end{equation*}
where
$V(\mathscr C_z^{\zb\xi}) $
denotes the $(d-2)$-dimensional volume of $\mathscr C_z^{\zb\xi}$.

The spherical transform $\mathcal U_z$ on the two-dimensional sphere $\S^2$ was first investigated  by Salman \cite{Sal15} in 2016.
He showed the injectivity of this transform $\mathcal U_z$ for smooth functions $f\in C^\infty(\S^2)$ that are supported inside the spherical cap $\{\zb\xi\in\S^2 \mid \xi_3 < z \}$ as well as an inversion formula.
This result was extended to the $d$-dimensional sphere in \cite{Sal16}, where also the smoothness requirement was lowered to $f\in C^1(\S^d)$.
A different approach was taken in \cite{Qu17}, where a relation with the classical Funk--Radon transform~$\mathcal U_0$ was established and also used for a characterization of the nullspace and range of $\mathcal U_z$ for $z<1$.

In the present paper, we extend the approach of \cite{Qu17} to the $d$-dimensional case.
We use a similar change of variables to connect the spherical transform $\mathcal U_z$ with the Funk--Radon transform $\mathcal F = \mathcal U_0$ in order to characterize its nullspace.
However, the description of the range requires some more effort, because the operator $\mathcal U_z$ is smoothing of degree ${({d-2})/{2}}$.

In the case $z=0$, we have sections of the sphere with hyperplanes through the origin, which are also called maximal subspheres of $\S^{d-1}$ or great circles on $\S^2$.
We obtain the classical Funk--Radon transform $\mathcal F = \mathcal U_0$.
The case $z=1$ corresponds to the subspheres containing the north pole $(0,\dots,0,1)^\top$.
This case is known as the spherical slice transform $\mathcal U_1$ and has been investigated since the early 1990s in \cite{AbDa93,Dah01,Hel11}.
However, unlike $\mathcal U_z$ for $z<1$, the spherical slice transform $\mathcal U_1$ is injective for $f\in L^\infty(\S^{d-1})$, see \cite{Rub15II}.
The main tool to derive this injectivity result of $\mathcal U_1$ is the stereographic projection, which turns the subspheres of $\S^{d-1}$ through the north pole into hyperplanes in $\R^{d-1}$ and thus connecting the spherical slice transform $\mathcal U_1$ to the Radon transform on $\R^{d-1}$.
For $z\to\infty$, we obtain vertical slices of the sphere, i.e., sections of the sphere with hyperplanes that are parallel to the north--south axis. 
This case is well-known for $\S^2$, see \cite{GiReSh94,ZaSc10,HiQu15circav,Rub18}.
In 2016, Palamodov \cite[Section~5.2]{Pal16} published an inversion formula for a certain nongeodesic Funk transform, which considers sections of the sphere with hyperplanes that have a fixed distance to the point $(0,\dots,0,z)^\top$.
This contains both the transform $\mathcal U_z$ as well as for $z=0$ the sections with subspheres having constant radius previously considered by Schneider \cite{Sch69} in 1969.

A key tool for analyzing the stability of inverse problems are Sobolev spaces, cf.\ \cite{JoLo01} for the Radon transform and \cite{HiQu15} for the Funk--Radon transform.
The Sobolev space $H^{s}(\S^{d-1})$ can be imagined as the space of functions $f\colon \S^{d-1}\to\C$ whose derivatives up to order $s$ are square-integrable, and we denote by $H^{s}_{\textrm{even}}(\S^{d-1})$ its restriction to even functions $f(\zb\xi)=f(-\zb\xi)$.
A thorough definition of $H^s(\S^{d-1})$ is given in \prettyref{sec:Sobolev-spaces}.
The behavior of the Funk--Radon transform was investigated by Strichartz \cite{Str81}, who found that the Funk--Radon transform $\mathcal F\colon H^s_{\textrm{even}}(\S^{d-1})\to H^{s+({d-2})/{2}}_{\textrm{even}}(\S^{d-1})$ is continuous and bijective.
In \prettyref{thm:Uz-Sobolev}, we show that basically the same holds for the generalized transform $\mathcal U_z$, i.e., the spherical transform $\mathcal U_z\colon H^s_{z}(\S^{d-1})\to H^{s+({d-2})/{2}}_{\textrm{even}}(\S^{d-1})$ is bounded and bijective and its inverse is also bounded, where $H^s_{z}(\S^{d-1})$ contains functions that are symmetric with respect to the point reflection in $(0,\dots,0,z)^\top$.

This paper is structured as follows.
In \prettyref{sec:def}, we give a short introduction to smooth manifolds and review the required notation on the sphere.
In \prettyref{sec:factorization}, we show the relation of $\mathcal U_z$ with the classical Funk--Radon transform $\mathcal F = \mathcal U_0$.
With the help of that factorization, we characterize the nullspace of the spherical transform $\mathcal U_z$ in \prettyref{sec:nullspace}.
In \prettyref{sec:function-spaces}, we consider Sobolev spaces on the sphere and show the continuity of certain multiplication and composition operators.
Then, in \prettyref{sec:range}, we prove the continuity of the spherical transform $\mathcal U_z$ in Sobolev spaces.
Finally, we include a geometric interpetation of our factorization result in \prettyref{sec:geometric-interpratation}.

\section{Preliminaries and definitions}
\label{sec:def}

\subsection{Manifolds}

We give a brief introduction to smooth manifolds, more can be found in \cite{AgFr02}.
We denote by $\R$ and $\C$ the real and complex numbers, respectively.
We define the $d$-dimensional Euclidean space $\R^{d}$
equipped with the scalar product $\left<\zb \xi,\zb \eta\right> = \sum_{i=1}^{d}\xi_i\eta_i$ and the norm $\norm{\zb\xi}=\left<\zb \xi,\zb \xi\right>^{1/2}$. 
We denote the unit vectors in $\R^{d}$ with $\zb \epsilon^i$, i.e.\ $\epsilon^{i}_j = \delta_{i,j}$, where $\delta$ is the Kronecker delta.
Every vector $\zb \xi\in\R^{d}$ can be written as $\zb \xi=\sum_{i=1}^{d}\xi_i \zb \epsilon^i$.

A diffeomorphism is a bijective, smooth mapping $\R^n\to \R^n$ whose inverse is also smooth.
We say a function is smooth if it has derivatives of arbitrary order.
An $n$-dimensional smooth manifold $ M$ without boundary is a subset of $\R^{d}$ such that for every $\zb \xi\in M$ there exists an open neighborhood $N(\zb \xi)\subset \R^{d}$ containing $\zb \xi$, an open set $U\subset\R^{d}$, and a diffeomorphism $m\colon U \to N(\zb \xi)$ such that 
\[
m(U\cap (\R^n\times \{0\}^{d-n})) = M \cap N(\zb \xi).
\]
We call $m\colon V\to M$ a map of the manifold $M$, where $V=U\cap (\R^n\times \{0\}^{d-n})$.
An atlas of $M$ is a finite family of maps $m_i\colon V_i\to M$, $i=1,\dots,l$, such that the sets $m_i (V_i)$  cover $M$.
We define the tangent space $T_{\zb\xi} M$ of $M$ at $\zb\xi\in M$ as the set of vectors $\zb x\in\R^{d}$ for which there exists a smooth path $\gamma\colon [0,1]\to M$ satisfying $\gamma(0) = \zb \xi$ and $\gamma'(0)=\zb x$.

A $k$-form $\omega$ on $M$ is a family $(\omega_{\zb\xi})_{\zb\xi \in M}$ of antisymmetric $k$-linear functionals 
$$\omega_{\zb\xi}\colon (T_{\zb\xi} M)^k \to \R,$$ 
where $(T_{\zb\xi} M)^k = T_{\zb\xi} M \times \dots \times T_{\zb\xi} M$.
Let $f\colon M \to N$ be a smooth mapping between the manifolds $M$ and $N$
and let $\omega$ be a $k$-form on $N$.
The pullback of $\omega$ is the $k$-form $f^*(\omega)$ on $M$ that is defined for any $\zb v^1,\dots,\zb v^k \in T_{\zb\xi} M$ by
\begin{equation*}
f^*(\omega) ([\zb v^i]_{i=1}^k)
= \omega \left(\left[\left.\frac{\dx}{\dx t}f\circ \gamma_i(t)\right|_{t=0}\right]_{i=1}^k\right),
\end{equation*}
where $\gamma_i\colon [0,1]\to M$ are smooth paths satisfying $\gamma_i(0) = \zb \xi$ and $\gamma_i'(0)=\zb v^i$.
If the smooth function $f\colon \R^{d}\to\R^{d}$ extends to the surrounding space, the pullback can be expressed as
\begin{equation}
f^*(\omega) ([\zb v^i]_{i=1}^k)
= \omega ([J_f\, \zb v^i]_{i=1}^k),
\label{eq:pullback-def}
\end{equation}
where $J_f$ denotes the Jacobian matrix of $f$.

An atlas $\{m_i\colon V_i \to M\}_{i=1}^l$ is called orientation of $M$ if for each $i,j$ the determinant of the Jacobian of $m_i^{-1} \circ m_j$ is positive wherever it exists.
A basis $\zb e^1,\dots,\zb e^n$ of the tangent space $T_{\zb\xi} M$ is oriented positively if for some $i$ with $\zb\xi \in m_i(V_i)$ the determinant of the matrix $[J_{m_i^{-1}}\, \zb e^1,\dots,J_{m_i^{-1}}\, \zb e^n]$ is positive.
Up to the multiplication of a constant real number depending only on ${\zb\xi}$, there exists only one $n$-form on an $n$-dimensional manifold.
Let $\zb e^1,\dots,\zb e^n$ be a positively oriented, orthonormal basis of the tangent space $T_{\zb\xi} M$.
Then the volume form $\dx M$ is the unique $n$-form on $M$ satisfying $\dx M(\zb e^1,\dots,\zb e^n) = 1$.

A set $\{\varphi_i\}_{i=1}^l$ of functions $\varphi_i\in C^\infty(M)$ is a partition of unity of the manifold $M$ with respect to an atlas $\{m_i\colon V_i \to M\}_{i=1}^l$ if
$\operatorname{supp}(\varphi_i)\subset m_i(V_i)$ for all $i$ and $\sum_{i=1}^{l} \varphi_i \equiv 1$ on $M$.
Then the integral of a $k$-form $\omega$ on $M$ is defined as 
\begin{equation*}
\int_M \omega
= \sum_{i=1}^{l} \int_{V_i} m_i^*(\varphi_i\cdot \omega)
= \sum_{i=1}^{l} \int_{V_i} c_i \d \R^n,
\end{equation*}
where the latter is the standard volume integral $\dx \R^n$ on $V_i\subset \R^n$ and the functions $c_i\colon V_i\to\R$ are uniquely determined by the condition $m_i^*(\varphi_i\cdot \omega) = c_i \d \R^n$.

Let $f\colon M\to N=f(M)$ be a diffeomorphism between the $n$-dimensional manifolds $M$ and $f(M)$, and let $\omega$ be an $n$-form on $N$. Then the 
substitution rule \cite[p.~94]{Jan01} holds,
	\begin{equation}
	\int_{f(M)} \omega
	= \int_M f^*(\omega).
	\label{eq:substitution-rule}	
	\end{equation}

\subsection{The sphere}
Let $d\ge3$.
The $(d-1)$-dimensional unit sphere 
\[\S^{d-1}=\{\zb\xi\in\R^{d}\mid\norm{\zb\xi}=1\}\] 
is a $(d-1)$-dimensional manifold in $\R^{d}$ with tangent space
\begin{equation*}
T_{\zb\xi} \S^{d-1} = \{\zb x\in\R^{d}\mid \left<\zb\xi,\zb x\right>=0\}
\end{equation*}
for $\zb\xi\in\S^{d-1}$.
We define an orientation on $\S^{d-1}$ by saying that a basis $[\zb x^i]_{i=1}^{d-1}$ of $T_{\zb\xi} \S^{d-1}$ is oriented positively if the determinant $\det [\zb \xi, \zb x^1, \dots, \zb x^{d-1}]>0$.
We denote the volume of the $(d-1)$-dimensional unit sphere $\S^{d-1}$ with
\begin{equation*}
\abs{\S^{d-1}}
= \int_{\S^{d-1}} \d\S^{d-1}
= \frac{2\pi^{d/2}}{\Gamma(d/2)}.
\end{equation*}

\subsubsection{The spherical transform}

Every $(d-2)$-dimensional subsphere of the sphere $\S^{d-1}$ is the intersection of $\S^{d-1}$ with a hyperplane, i.e. 
\begin{equation*}
S(\zb{\xi},t) = \left\{
\boldsymbol{\eta}\in\S^{d-1} \mid \left\langle
\boldsymbol{\xi},\boldsymbol{\eta}\right\rangle =t\right\},
\end{equation*} 
where 
$\zb{\xi}\in\S^{d-1}$ is the normal vector of the hyperplane 
and $t\in[-1,1]$ is the signed distance of the hyperplane to the origin.
We define an orientation on the subsphere $S(\zb{\xi},t)$ by saying that a basis $[\zb e^i]_{i=1}^{d-2}$ of the tangent space $T_{\zb \eta} S(\zb{\xi},t)$  is oriented positively if
\[
\det \left(\zb\eta, \zb\xi, \zb e^1,\dots,\zb e^{d-2} \right) >0.
\]

In the following, we consider subspheres of $\S^{d-1}$ whose hyperplanes have a common point 
located in the interior of the unit ball. 
Because of the rotational symmetry, we can assume that this point lies on the positive $\xi_{d}$ axis.
For $z\in[0,1)$, we consider the point
\[
z\zb\epsilon^{d} = (0,\dots,0,z)^\top.
\] 
The $(d-2)$-dimensional subsphere that is located in one hyperplane together with $z\zb\epsilon^{d}$ can be described by $S(\zb \xi,t)$ with $\zb \xi\in\S^{d-1}$ and $t = \left<\zb \xi, z\zb\epsilon^{d}\right> = z\xi_{d}$.
We define the subsphere
\begin{equation*}
\mathscr C_z^{\zb\xi}
= S (\zb \xi, z\xi_{d}).
\end{equation*}
The $(d-2)$-dimensional sphere $\mathscr C_z^{\zb\xi}$ has radius $\sqrt{1-z^2\xi_{d}^2}$ and volume
\begin{equation*}
V(\mathscr C_z^{\zb\xi}) =
\abs{\S^{d-2}} \left({1-z^2 \xi_{d}^2}\right)^{(d-2)/2}.
\end{equation*}
For a continuous function $f\colon\S^{d-1}\to\C$, we define its {spherical transform} $\mathcal U_z f$ by
\begin{equation}
\mathcal U_z f (\zb{\xi}) 
= \frac{1}{V(\mathscr C_z^{\zb\xi})} 
\int_{\mathscr C_z^{\zb\xi}} f \d \mathscr C_z^{\zb\xi}
,\qquad \zb{\xi} \in\S^{d-1},
\label{eq:spherical-transform-def}
\end{equation}
which computes the mean values of $f$ along the subspheres $\mathscr C_z^{\zb\xi}$.

\begin{Remark}
	The subspheres $\mathscr C_z^{\zb\xi}$, along which we integrate, can also be imagined in the following way.
	The centers of the subspheres $\mathscr C_z^{\zb\xi}$ are located on a sphere 
    that contains the origin and $z\zb\epsilon^{d}$ and is rotationally symmetric about the North--South axis.
	This can be seen as follows.
	The center of the sphere $\mathscr C_z^{\zb\xi}$ is given by $z\xi_{d}\zb\xi$.
	Then the distance of $z\xi_{d}\zb\xi$ to the point $\frac{z}{2}\,\zb\epsilon^{d}$ reads
	\begin{align*}
	\left\lVert{z\xi_{d}\zb\xi - \frac{z}{2}\,\zb\epsilon^{d}}\right\rVert^2
	&= \sum_{i=1}^{d-1} \xi_i^2 z^2 \xi_{d}^2 + \left(z\xi_{d}^2-\frac{z}{2}\right)^2
	\\
	&= (1-\xi_{d}^2) z^2 \xi_{d}^2 + \left(z\xi_{d}^2-\frac{z}{2}\right)^2
	\\&=\left(\frac{z}{2}\right)^2,
	\end{align*}
	which is independent of $\zb\xi$.
	So, the centers of $\mathscr C_z^{\zb\xi}$ are located on a sphere with center $\frac{z}{2}\zb\epsilon^{d}$ and radius $\frac{z}{2}$.
\end{Remark}

\subsubsection{The Funk--Radon transform}
\label{sec:FRT}

Setting the parameter $z=0$,
the point $0\zb\epsilon^d=(0,\dots,0)^\top$ is the center of the sphere $\S^{d-1}$.
Hence, the spherical transform $\mathcal U_0$ integrates along all maximal subspheres of the sphere $\S^{d-1}$.
This special case is called the {Funk--Radon transform}
\begin{equation}
\mathcal Ff (\zb \xi)
= \mathcal U_0 f (\zb \xi)
= \frac{1}{\abs{\S^{d-2}}} \int_{\mathscr C_0^{\zb\xi}} f \d \mathscr C_0^{\zb\xi}, \qquad
\zb{\xi}\in\S^{d-1},
\label{eq:FRT}
\end{equation} 
which is also known by the terms {Funk transform}, {Minkowski--Funk transform} or {spherical Radon transform}, where the latter term occasionally also refers to means over $({d-1})$-dimensional spheres in $\R^{d}$, cf.\ \cite{Qui82}.

\section{Relation with the Funk--Radon transform}
\label{sec:factorization}

In this section, we 
show a connection between the spherical transform $\mathcal U_z$ and the Funk--Radon transform $\mathcal F$.

\subsection{Two mappings on the sphere}
\label{sec:hg}

Let $z\in(-1,1)$.
We define the transformations $\zb h_z, \zb g_z\colon \S^{d-1}\to \S^{d-1}$ by
\begin{equation}
\zb h_z (\zb \eta)
= \sum_{i=1}^{d-1} \frac{\sqrt{1-z^2}}{1+z\eta_{d}} \eta_i \zb \epsilon^i + \frac{z+\eta_{d}}{1+z\eta_{d}} \zb \epsilon^{d}
\label{eq:h_z}
\end{equation}
and
\begin{equation}
\zb g_z (\zb \xi)
= \frac{1}{\sqrt{1-z^2\xi_{d}^2}} 
\left(\sum_{i=1}^{d-1} \xi_i \zb \epsilon^i + \sqrt{1-z^2} \xi_{d} \zb \epsilon^{d}\right).
\label{eq:g_z}
\end{equation}

\begin{Remark}
\label{rem:gh-injective}
The definitions of both $\zb h_z$ and $\zb g_z$ rely only on the $d$-th coordinate.
The values in the other coordinates are just multiplied with the same factor in order to make the vectors stay on the sphere.
Furthermore, the transformations $\zb h_z$ and $\zb g_z$ are bijective with their respective inverses given by
\begin{equation}
\zb h_z^{-1}({\zb\omega})
= \zb h_{-z} (\zb\omega)
=\sum_{i=1}^{d-1} \frac{\sqrt{1-z^2}}{1-z\omega_{d}}\omega_i \zb\epsilon^i +
\frac{\omega_{d}-z}{1-z\omega_{d}}\zb\epsilon^{d}
\label{eq:h^-1}
\end{equation}
and
\begin{equation}
\zb g_z^{-1}(\zb\omega)
=\zb g_{\frac{\mathrm iz}{\sqrt{1-z^2}}}(\zb\omega)
=\frac{1}{\sqrt{1-z^2+z^2\omega_{d}^2}} \left(\sum_{i=1}^{d-1} \sqrt{1-z^2}\,\omega_i \zb\epsilon^i
+\omega_{d}\zb\epsilon^{d}\right).
\label{eq:g^-1}
\end{equation}
The computation of the inverses is straightforward and therefore omitted here.
\end{Remark}

The following lemma shows that the inverse of $\zb h_z$ applied to the subsphere $\mathscr C_z^{\zb\xi}$ yields a maximal subsphere of $\S^{d-1}$ with normal vector $\zb g_z(\zb\xi)$.

\begin{Lemma}
	Let $z\in (-1,1)$ and $\zb\xi\in \S^{d-1}$.
	Then 
	\begin{equation}
	\zb h_z^{-1} (\mathscr C_z^{\zb\xi}) = \mathscr C_0^{\zb g_z(\zb\xi)}.
	\label{eq:CD}
	\end{equation}
\end{Lemma}
\begin{proof}
	Let $\zb\eta\in\S^{d-1}$.
	Then $\zb \eta$ lies in $\zb h_z^{-1} (\mathscr C_z^{\zb\xi})$ if and only if $\zb h_z(\zb \eta) \in \mathscr C_z^{\zb\xi}$, i.e.,
	\begin{equation*}
	\left\langle\zb h_z(\zb{\eta}), \zb{\xi}\right\rangle =z\xi_{d}.
	\end{equation*}
	By the definition of $\zb h_z$ in \eqref{eq:h_z}, we have
	\begin{equation*}
	\sum_{i=1}^{d-1} \frac{\sqrt{1-z^2}}{1+z\eta_{d}} \eta_i \xi_i + \frac{z+\eta_{d}}{1+z\eta_{d}} \xi_{d}
	= z\xi_{d}.
	\end{equation*}
	After subtracting the right-hand side from the last equation, we have
	\begin{equation*}
	\sum_{i=1}^{d-1} \frac{\sqrt{1-z^2}}{1+z\eta_{d}} \eta_i \xi_i + \frac{1-z^2}{1+z\eta_{d}} \eta_{d}\xi_{d}
	= 0.
	\end{equation*}
	Multiplication with $(1+z\eta_{d})(1-z^2)^{-1/2}(1-z^2\xi_{d}^2)^{-1/2}$ yields
	\begin{equation*}
	\sum_{i=1}^{d-1} \frac{1}{\sqrt{1-z^2\xi_{d}^2}} \eta_i \xi_i + \frac{\sqrt{1-z^2}}{\sqrt{1-z^2\xi_{d}^2}} \eta_{d}\xi_{d}
	= 0,
	\end{equation*}
	which is equivalent to $\langle \zb \eta, \zb g_z(\zb \xi)\rangle =0$,
	so we obtain that $\zb \eta \in \mathscr C_0^{\zb g_z(\zb\xi)}$.
\end{proof}

\begin{Lemma}
	\label{lem:pullback}
    Let $z\in (-1,1)$ and $\zb\xi\in \S^{d-1}$.
	Denote by $\dx \mathscr C_z^{\zb\xi}$ and $\dx \mathscr C_0^{\zb g_z(\zb\xi)}$ the volume forms on the manifolds $ \mathscr C_z^{\zb\xi}$ and $\mathscr C_0^{\zb g_z(\zb\xi)}$, respectively.
	Then the following relation between the pullback of the volume form $\dx \mathscr C_z^{\zb\xi}$ over $\zb h_z$ and $\dx \mathscr C_0^{\zb g_z(\zb\xi)}$ holds.
	For $\zb \eta \in \mathscr C_0^{\zb g_z(\zb\xi)}$, we have
	\begin{equation}
	\zb h_z^*(\dx \mathscr C_z^{\zb\xi})
	= \left(\frac{\sqrt{1-z^2}}{1+z\eta_{d}}\right)^{d-2} \d \mathscr C_0^{g(\zb\xi)}.
	\label{eq:pullback}
	\end{equation}
	Furthermore, we have for the volume form $\d \S^{d-1}$ on the sphere
	\begin{equation}
	\zb h_z^*(\dx \S^{d-1})
	= \left(\frac{\sqrt{1-z^2}}{1+z\eta_{d}}\right)^{d-1} \d \S^{d-1}.
	\label{eq:pullback-sd}
	\end{equation}
\end{Lemma}
\begin{proof}
	We compute the Jacobian $J_{\zb h_z}$ of $\zb h_z$, which comprises the partial derivatives of $\zb h_z$.
	For all $l,m \in \{1,\dots,d-1\}$, we have
	\begin{equation}
	\begin{aligned}
	\frac{\partial [\zb h_z]_l}{\partial \eta_m}
	&= \frac{\sqrt{1-z^2}}{1+z\eta_{d}}\delta_{l,m},\qquad\qquad\qquad
	&
	\frac{\partial [\zb h_z]_l}{\partial \eta_{d}}
	&= -\eta_l\frac{z\sqrt{1-z^2}}{(1+z\eta_{d})^2},
	\\
	\frac{\partial [\zb h_z]_{d}}{\partial \eta_m}
	&= 0,
	&
	\frac{\partial [\zb h_z]_{d}}{\partial \eta_{d}}
	&= \frac{1-z^2}{(1+z\eta_{d})^2}.
	\end{aligned}
	\label{eq:Jacobian-h}
	\end{equation}
	Let $[\zb e^i]_{i=1}^{d-2}$ be an orthonormal basis of the tangent space $T_{\zb\eta} \mathscr C_0^{\zb g_z(\zb\xi)}$.
	Then $J_{\zb h_z} \zb e^i\in T_{\zb h_z(\zb \eta)} \mathscr C_z^{\zb\xi}$ for $i=1,\dots,d-1$ is given by
	\begin{align*}
	J_{\zb h_z}\zb e^i
	&= \sum_{l=1}^{d-1} \left(\frac{\sqrt{1-z^2}}{1+z\eta_{d}} e^i_l -\eta_l\frac{z\sqrt{1-z^2}}{(1+z\eta_{d})^2} e^i_{d} \right) \zb\epsilon^{l} 
	+ \frac{1-z^2}{(1+z\eta_{d})^2} e^i_{d} \zb\epsilon^{d} .
	\end{align*}
	Hence, we have for all $i,j \in \{1,\dots,d-2\}$
	\begin{align*}
	\left< J_{h_z}\zb e^i,J_{h_z}\zb e^j\right>
	={}& \sum_{l=1}^{d-1} \left(\frac{1-z^2}{(1+z\eta_{d})^2} e^i_le^j_l -\frac{z(1-z^2)}{(1+z\eta_{d})^3} \eta_l \left(e^i_le^j_{d}+e^j_le^i_{d}\right)\right.
	\\&\hspace{3em}
	\left. + \frac{z^2(1-z^2)}{(1+z\eta_{d})^4} \eta_l^2 e^i_{d}e^j_{d} \right) 
	+ \frac{(1-z^2)^2}{(1+z\eta_{d})^4} e^i_{d} e^j_{d}.
	\end{align*}
	Expanding the sum, we obtain
	\begin{align*}
	\left< J_{\zb h_z}\zb e^i,J_{\zb h_z}\zb e^j\right>
	={}& \frac{1-z^2}{(1+z\eta_{d})^2} \sum_{l=1}^{d-1}e^i_le^j_l
	-\frac{z(1-z^2)}{(1+z\eta_{d})^3}  
	\left(e^j_{d} \sum_{l=1}^{d-1}\eta_l e^i_l + e^i_{d}\sum_{l=1}^{d-1} \eta_l e^j_l\right)
	\\&
	+ \frac{z^2(1-z^2)}{(1+z\eta_{d})^4} e^i_{d}e^j_{d} \sum_{l=1}^{d-1} \eta_l^2 
	+ \frac{(1-z^2)^2}{(1+z\eta_{d})^4} e^i_{d} e^j_{d}.
	\end{align*}
	Since the vectors $\zb e^i$ and $\zb e^j$ are elements of an orthonormal basis, we have $\left<\zb e^i,\zb e^j\right> = \sum_{l=1}^{d} e^i_l e^j_l = \delta_{i,j}$. 
	Furthermore, we know that $\left<\zb e^i,\zb\eta\right>=\left<\zb e^j,\zb\eta\right> = 0$ because $\zb e^i$ and $\zb e^j$ are in the tangent space $T_{\zb\eta} \mathscr C_0^{g(\zb\xi)} \subset T_{\zb\eta} \S^{d-1}$,
	and also $\norm{\zb \eta}^2=1$.
	Hence, we have
	\begin{align*}
	&\left< J_{\zb h_z}\zb e^i,J_{\zb h_z}\zb e^j\right>
	\\
	={}& \frac{1-z^2}{(1+z\eta_{d})^2} (\delta_{i,j} - e^i_{d}e^j_{d})
	+2\frac{z(1-z^2)}{(1+z\eta_{d})^3} \eta_{d}e^i_{d}e^j_{d}
	\\&
	+ \frac{z^2(1-z^2)}{(1+z\eta_{d})^4} e^i_{d}e^j_{d} (1-\eta_{d}^2)
	+ \frac{(1-z^2)^2}{(1+z\eta_{d})^4} e^i_{d} e^j_{d}
	\\
	={}& \frac{1-z^2}{(1+z\eta_{d})^4} e^i_{d}e^j_{d} 
	\left(-(1+z\eta_{d})^2 +2z\eta_{d}(1+z\eta_{d})
	+ z^2 (1-\eta_{d}^2) +1-z^2
	\right) 
	\\&
	+ \frac{1-z^2}{(1+z\eta_{d})^2} \delta_{i,j}
	\\
	={}& \frac{1-z^2}{(1+z\eta_{d})^2} \delta_{i,j}.
	\end{align*}
	The above computation shows that the vectors $\{J_{\zb h_z} \zb e^i\}_{i=1}^{d}$ are orthogonal with length $\norm{J_{\zb h_z} \zb e^i} = {\sqrt{1-z^2}}/({1+z\eta_{d}})$.
	By the definition of the pullback in \eqref{eq:pullback-def} and the fact that the volume form $\dx \mathscr C_z^{\zb\xi}$ is a multilinear $(d-2)$-form, we obtain
	\begin{equation}
	\zb h_z^*(\dx \mathscr C_z^{\zb\xi}) ([\zb e^i]_{i=1}^{d-2})
	= \dx \mathscr C_z^{\zb\xi} ([J_{\zb h_z} \zb e^i]_{i=1}^{d-2})
	= \left(\frac{\sqrt{1-z^2}}{1+z\eta_{d}}\right)^{d-2}.
    \label{eq:pullback-1}
	\end{equation}
	If we set $[\zb e^i]_{i=1}^{d-1}$ as a basis of the tangent space $T_{\zb\eta}(\S^{d-1})$ in order to obtain \eqref{eq:pullback-sd}, the previous calculations still hold except that the exponent $d-2$ is replaced by $d-1$ in equation \eqref{eq:pullback-1}.
	
	Finally, we prove that the basis $\left[J_{\zb h_z}\zb e^1,\dots, J_{\zb h_z}\zb e^{d} \right]$ of $T_{\zb h_z(\zb\eta)}\mathscr C^{\zb\xi}_z$ is oriented positively, i.e., that
	\[
	d(z):=\det\left(\zb h_z(\zb\eta), \zb\xi, J_{\zb h_z}\zb e^1,\dots, J_{\zb h_z}\zb e^{d} \right)>0.
	\]
	By the formula \eqref{eq:Jacobian-h} of $J_{\zb h_z}$, the function $d\colon [0,1)\to\R$
	is continuous, 
	and it satisfies
	\[
	d(0)
	=\det\left( \zb h_0(\zb\eta), \zb\xi, J_{\zb h_0}\zb e^1,\dots, J_{\zb h_0}\zb e^{d-2} \right)
	=\det\left( \zb\eta, \zb g_0(\zb\xi), \zb e^1,\dots, \zb e^{d-2} \right) >0
	\]
	since both $\zb h_0$ and $\zb g_0$ are equal to the identity map and we assumed the orthonormal basis $\left[\zb e^1,\dots, \zb e^{d-2} \right]$ be oriented positively.
	By the orthogonality of the vectors $\zb\xi$, $\zb h_z(\zb\eta)$ and $J_{\zb h_z}\zb e^i$,
	we see that $d(z)$ vanishes nowhere and,
	hence, we obtain that $d(z)>0$ for all $z\in[0,1)$.
	The assertion follows by the uniqueness of the volume form $\dx \mathscr C_0^{\zb g_z(\zb\xi)}$.
\end{proof}

\subsection{Factorization}

Let $z\in(-1,1)$ and $f\in C(\S^{d-1})$. We define the two transformations $\mathcal M_z,\, \mathcal N_z\colon C(\S^{d-1})\to C(\S^{d-1})$ by
\begin{equation}
\mathcal M_zf (\zb{\xi}) =
\left(\frac{\sqrt{1-z^2}}{1+z\xi_{d}}\right)^{d-2} f \circ \zb h_z (\zb\xi) 
,\qquad \zb{\xi}\in\S^{d-1}
\label{eq:M}
\end{equation}
and
\begin{equation}
\mathcal N_zf (\zb{\xi}) =
(1-z^2\xi_{d}^2)^{-\frac{d-2}{2}}\,
f \circ \zb g_z ( \zb \xi)
,\qquad \zb{\xi}\in\S^{d-1}.
\label{eq:N}
\end{equation}

\begin{Remark}
	The transformations $\mathcal M_z$ and $\mathcal N_z$ are inverted by
	\begin{equation}
	f(\zb\eta)
	= \left(\frac{\sqrt{1-z^2}}{1-z\eta_{d}}\right)^{d-2} \mathcal M_z f(\zb h_z^{-1}(\zb\eta))
	= \mathcal M_{-z} \mathcal M_z f(\zb\eta),
	\label{eq:M-1}
	\end{equation}
	and
	\begin{equation}
	f(\zb\eta)
	= \left(\frac{1-z^2}{1-(1-z^2)\eta_{d}^2}\right)^{\frac{d-2}{2}} \mathcal N_z (\zb g_z^{-1}(\zb\eta)),
	\qquad \zb\eta\in\S^{d-1},
	\label{eq:N-1}
	\end{equation}
	respectively.
\end{Remark}

Now we are able to prove our main theorem about the factorization of the spherical transform $\mathcal U_z$.

\begin{Theorem}
	Let $z\in[0,1)$.
	Then the factorization of the spherical transform
	\begin{equation}
	\mathcal U_z 
	= \mathcal N_z \mathcal F \mathcal M_z
	\label{eq:spherical-transform-decomposition}
	\end{equation}
	holds, where $\mathcal F$ is the Funk--Radon transform \eqref{eq:FRT}.
	\label{thm:spherical-transform}
\end{Theorem}

\begin{proof}
	Let $f\in C(\S^2)$ and $\zb \xi\in\S^{d-1}$.
	By the definition of $\mathcal U_z$ in \prettyref{eq:spherical-transform-def}, we have
	\begin{equation}
	\abs{\S^{d-2}} (1-z^2\xi_{d}^2)^{\frac{d-2}{2}}\, \mathcal U_z f(\zb{\xi})
	=\int_{ \mathscr C_z^{\zb\xi}} f \d \mathscr C_z^{\zb\xi}.
	\label{eq:spherical-transform-1}
	\end{equation}
	Then we have by the substitution rule \eqref{eq:substitution-rule} 
	\begin{equation*}
	\int_{\mathscr C_z^{\zb\xi}} f \d \mathscr C_z^{\zb\xi}
	= \int_{\zb h_z^{-1}(\mathscr C_z^{\zb\xi})} (f\circ\zb h_z) \, \zb h_z^*(\dx \mathscr C_z^{\zb\xi}).
	\end{equation*}
	By \eqref{eq:pullback} and \eqref{eq:CD}, we obtain
	\begin{equation*}
	\int_{\mathscr C_z^{\zb\xi}} f \d \mathscr C_z^{\zb\xi}
	= \int_{\mathscr C_0^{\zb g_z(\zb\xi)}} f(\zb h_z(\zb \eta)) \left(\frac{\sqrt{1-z^2}}{1+z\eta_{d}}\right)^{d-2} \d \mathscr C_0^{\zb g_z(\zb\xi)} (\zb\eta).
    \end{equation*}
    By the definition of $\mathcal M_z$ in \eqref{eq:M}, we see that
    \begin{equation*}
    \int_{\mathscr C_z^{\zb\xi}} f \d \mathscr C_z^{\zb\xi}
	= \int_{\mathscr C_0^{\zb g_z(\zb\xi)}} \mathcal M_z f \d \mathscr C_0^{\zb g_z(\zb\xi)}.
	\end{equation*}
	The definition of the Funk--Radon transform \eqref{eq:FRT} shows that
	\begin{equation*}
	\int_{\mathscr C_z^{\zb\xi}} f \d \mathscr C_z^{\zb\xi}
	= \mathcal F \mathcal M_z f(\zb g_z(\zb \xi)),
	\end{equation*}
	which implies \eqref{eq:spherical-transform-decomposition}.
\end{proof}

The factorization theorem \ref{thm:spherical-transform} enables us to investigate the properties of the spherical transform $\mathcal U_z$.
Because the operators $\mathcal M_z$ and $\mathcal N_z$ are relatively simple, we can transfer many properties from the Funk--Radon transform, which has been studied by many authors already, to the spherical transform $\mathcal U_z$.

\section{Nullspace}

\label{sec:nullspace}

With the help of the factorization \eqref{eq:spherical-transform-decomposition} obtained in the previous section, we obtain the following characterization of the nullspace of the spherical transform $\mathcal U_z$.

\begin{Theorem}
	\label{thm:nullspace}
	Let $z\in[0,1)$ and $f\in C(\S^{d-1})$.
	Then $\mathcal U_zf=0$ if and only if
	\begin{equation}
	f(\zb{\omega})
	=-\left(\frac{1-z^2}{1-2z\omega_{d}+z^2}\right)^{d-2} f \circ\zb r_z(\zb\omega)
	,\qquad \zb\omega\in\S^{d-1},
	\end{equation}
	where $\zb r_z\colon \S^{d-1}\to\S^{d-1}$ is given by
	\begin{equation}
	\label{eq:point-reflection}
	\zb r_z(\zb\omega)= \sum_{i=1}^{d-1} \frac{z^2-1}{1+z^2-2z\omega_{d}}\omega_i\zb\epsilon^i
	+ \frac{2z-z^2\omega_{d}-\omega_{d}}{1+z^2-2z\omega_{d}}\zb\epsilon^{d}.
	\end{equation}
\end{Theorem}

\begin{proof}
	Let $f\in C(\S^{d-1})$. Since the operator $\mathcal N_z$ is bijecive by \prettyref{rem:gh-injective}, we see that $\mathcal U_zf = \mathcal N_z \mathcal F \mathcal M_z f = 0$ if and only if $\mathcal F \mathcal M_z f=0$.
	The nullspace of the Funk--Radon transform $\mathcal F$ consists of the odd functions, cf.\ \cite[Proposition 3.4.12]{Gro96}, so we obtain
$$\mathcal M_z f(\zb \eta) = -\mathcal M_z f(-\zb \eta),\qquad \zb\eta\in\S^{d-1}.$$
By the definition of $\mathcal M_z$ in \eqref{eq:M}, we have
$$ \left(\frac{\sqrt{1-z^2}}{1+z\eta_{d}}\right)^{d-2} f \circ \zb h_z (\zb\eta) 
= -\left(\frac{\sqrt{1-z^2}}{1-z\eta_{d}}\right)^{d-2} f \circ \zb h_z (-\zb\eta). $$
We substitute $\zb\omega=\zb h_z(\zb\eta)$ and obtain
\begin{align*}
f (\zb\omega) 
&=- \left(\frac{1+z\frac{\omega_{d}-z}{1-z\omega_{d}}}{1-z\frac{\omega_{d}-z}{1-z\omega_{d}}}\right)^{d-2} f \circ \zb h_z (-\zb h_z^{-1}(\zb\omega))
\\&
=- \left(\frac{1-z\omega_{d}+z(\omega_{d}-z)}{1-z\omega_{d}-z{(\omega_{d}-z})}\right)^{d-2} f \circ \zb h_z (-\zb h_z^{-1}(\zb\omega))
\\&
=- \left(\frac{1-z^2}{1-2z\omega_{d}+z^2}\right)^{d-2} f \circ \zb h_z (-\zb h_z^{-1}(\zb\omega)).
\end{align*}
In order to show that $\zb r_z = \zb h_z(-\zb h_z^{-1})$,
we compute the $d$-th component
\begin{align*}
[\zb h_z(-\zb h_z^{-1}(\zb\omega))]_{d}
= \frac{z+\frac{z-\omega_{d}}{1-z\omega_{d}}}
{1+z\frac{z-\omega_{d}}{1-z\omega_{d}}}
= \frac{z-z^2\omega_{d}+{z-\omega_{d}}}
{1-z\omega_{d}+z^2-z\omega_{d}}
= \frac{2z-z^2\omega_{d}-\omega_{d}}
{1-2z\omega_{d}+z^2}.
\end{align*}
For $i\in\{1,\dots,d-1\}$, we have
\begin{equation*}
[\zb h_z(-\zb h_z^{-1}(\zb\omega))]_{i}
= \frac{\sqrt{1-z^2}}{1-z\frac{\omega_{d}-z}{1-z\omega_{d}}} \frac{-\sqrt{1-z^2}}{1-z\omega_{d}} \omega_i
= \frac{z^2-1}{1-2z\omega_{d}+z^2}\, \omega_i.\qedhere
\end{equation*}
\end{proof}

\begin{Remark}
	The map $\zb r_z$ from \eqref{eq:point-reflection}
	is the point reflection of the sphere $\S^{d-1}$ about the point $z\zb\epsilon^{d}$. 
	This can be seen as follows.
    Let $\zb\omega\in\S^{d-1}$.
	The vectors $\zb\omega-z\zb\epsilon^{d}$ and $\zb r_z(\zb\omega)-z\zb\epsilon^{d}$ are parallel if for all $i\in\{1,\dots,d\}$
	\begin{equation*}
	\frac{[\zb r_z(\zb\omega)]_i}{\omega_i}
	=\frac{[\zb r_z(\zb\omega)]_{d}-z}{\omega_{d}-z}.
	\end{equation*}
	We have
	\begin{align*}
	\frac{\omega_i}{[\zb r_z(\zb\omega)]_i} 
	\frac{[\zb r_z(\zb\omega)]_{d}-z}{\omega_{d}-z}
	&=\frac{2z-z^2\omega_{d}-\omega_{d} -z(1+z^2-2z\omega_{d})}{(z^2-1)(\omega_{d}-z)}
	\\&=\frac{z+z^2\omega_{d}-\omega_{d} -z^3}{(z^2-1)(\omega_{d}-z)}
	=1,
	\end{align*}
	provided all denominators are nonzero.
\end{Remark}

\section{Function spaces on the sphere}
\label{sec:function-spaces}

Before we can state the range of the spherical transform $\mathcal U_z$, we have to introduce some function spaces on the sphere $\S^{d-1}$.
The Hilbert space $L^2(\S^{d-1})$ comprises all square-integrable functions with the inner product of two functions $f,g\colon\S^{d-1}\to\C$
$$
\left< f,g \right>_{L^2(\S^{d-1})}
= \int_{\S^{d-1}} f(\zb\xi)\ \overline{g(\zb\xi)} \d \S^{d-1}(\zb\xi)
$$ 
and the norm $\norm{f}_{L^2(\S^{d-1})}={\left<f,f\right>_{L^2(\S^{d-1})}^{1/2}}$.

\subsection{The space $C^s(\S^{d-1})$ and differential operators on the sphere}

For brevity, we denote by $\partial_i = \frac{\partial}{\partial x_i}$ the partial derivative with respect to the $i$-th variable.
We extend a function $f\colon\S^{d-1}\to\C$ to the surrounding space $\R^{d} \setminus\{\zb 0\}$ by setting
\begin{equation*} 
f^\bullet(\zb x) = f\left(\frac{\zb x}{\norm{\zb x}}\right) 
,\qquad \zb x\in\R^{d} \setminus\{\zb 0\}.
\end{equation*}
The surface gradient $\zb\nabla^\bullet$ on the sphere is the orthogonal projection of the gradient $\zb\nabla=(\partial_1,\dots,\partial_{d})^\top$ onto the tangent space of the sphere.
For a differentiable function $f\colon\S^{d-1}\to\C$, we have
\begin{equation*}
\zb\nabla^\bullet f (\zb\xi) = \zb\nabla f^\bullet(\zb\xi)
,\qquad \zb\xi\in\S^{d-1}.
\end{equation*}
In a similar manner, the restriction of the Laplacian
\begin{equation*}
\Delta = \partial_1^2 + \dots + \partial_{d}^2
\end{equation*}
to the sphere is known as the Laplace--Beltrami operator \cite[(§14.20)]{Mue98}
\begin{equation*}
\Delta^\bullet f (\zb\xi) = \Delta f^\bullet(\zb\xi)
,\qquad \zb\xi\in\S^{d-1}.
\end{equation*}
For a multi-index $\zb\alpha=(\alpha_1,\dots,\alpha_{d}) \in\N^{d}$, we define its norm $\norm{\zb\alpha}_1 = \sum_{i=1}^{d}\abs{\alpha_i}$ and the differential operator $D^{\zb\alpha} = \partial_1^{\alpha_1} \cdots \partial_{d}^{\alpha_{d}}$.
Let $s\in\N$. 
We denote by $C^s(\S^{d-1})$ the space of functions $f\colon \S^{d-1}\to\C$ whose extension $f^\bullet$ has continuous derivatives up to the order $s$ with the norm
\begin{equation*}
\norm{f}_{C^{s}(\S^{d-1})}
= \max_{\norm{\zb\alpha}_1 \leq s}\, \sup_{\zb\xi\in\S^{d-1}} \abs{D^{\zb\alpha} f^\bullet(\zb\xi)}.
\end{equation*}
The space $C^0(\S^{d-1}) = C(\S^{d-1})$ is the space of continuous functions with the uniform norm.
The definition implies for $f\in C^{s+1}(\S^{d-1})$
\begin{equation}
\norm{f}_{C^s(\S^{d-1})} \leq \norm{f}_{C^{s+1}(\S^{d-1})}.
\label{eq:infty-norm-1}
\end{equation}
We define the space $C^s(\S^{d-1}\to\R^{d})$ of vector fields $\zb f\colon\S^{d-1}\to\R^{d}$ with the norm as the Euclidean norm over its component functions, i.e., for $\zb f(\zb\xi) =  [f_i(\zb\xi)]_{i=1}^{d}$ we set
\begin{equation*}
\norm{\zb f}_{C^s(\S^{d-1}\to\R^{d})}
= \sqrt{\sum_{i=1}^{d} \norm{f_i}_{C^s(\S^{d-1})}^2}.
\label{eq:infty-norm-vector}
\end{equation*}
We see that  for $f\in C^{s+1}(\S^{d-1})$
\begin{equation}
\norm{\zb\nabla^\bullet f}_{C^{s}(\S^{d-1}\to\R^d)}^2
= \sum_{i=1}^{d} \norm{\partial_i f^\bullet}_{C^{s}(\S^{d-1})}^2
\leq \sum_{i=1}^{d} \norm{ f}_{C^{s+1}(\S^{d-1})}^2 
= d\, \norm{f}_{C^{s+1}(\S^{d-1})}^2.
\label{eq:infty-norm-gradient}
\end{equation}

\subsection{Sobolev spaces}

\label{sec:Sobolev-spaces}

We give a short introduction to Sobolev spaces on the sphere based on \cite{AtHa12} (see also \cite{DaXu2013,Mic13}).
We define the Legendre polynomial $P_{n,d}$ of degree $n\in\N$ and in dimension $d$ by \cite[(2.70)]{AtHa12}
\begin{equation*}
P_{n,d}(t)
= (-1)^n\, \frac{(d-3)!!}{(2n+d-3)!!}\, (1-t^2)^{\frac{3-d}{2}} \left(\frac{\dx}{\dx t}\right)^n (1-t^2)^{n+\frac{d-3}{2}}
,\qquad t\in[-1,1].
\end{equation*}
For $f\in L^2(\S^{d-1})$ and $n\in\N$, we define the projection operator
$$
\mathcal P_{n,d}f(\zb\xi) = \frac{N_{n,d}}{\abs{\S^{d-1}}} 
\int_{\S^{d-1}} f(\zb \eta)\,  P_{n,d}(\left<\zb\xi,\zb \eta\right>) \d \S^d(\zb \eta)
,\qquad\zb\xi\in\S^{d-1},
$$
where
$$
N_{n,d} = \dim(\mathcal P_{n,d} (L^2(\S^{d-1}))) = \frac{(2n+d-2)(n+d-3)!}{n!(d-2)!}.
$$
Note that $\mathcal P_{n,d}$ is the $L^2(\S^{d-1})$-orthogonal projection onto the pairwise orthogonal spaces $
\mathcal P_{n,d} (L^2(\S^{d-1}))$ of harmonic polynomials that are homogeneous of degree $n$ restricted to the sphere $\S^{d-1}$.
Every function $f\in L^2(\S^{d-1})$ can be written as the Laplace series
$$
f = \sum_{n=0}^{\infty} \mathcal P_{n,d} f.
$$
We define the Sobolev space
$H^s(\S^{d-1})$ 
of smoothness $s\ge0$
as the space of all functions $f\in L^2(\S^{d-1})$ with finite Sobolev norm \cite[(3.98)]{AtHa12}
\begin{align}
\label{eq:sobolev-norm}
\norm{f}_{H^s(\S^{d-1})}
&= \sqrt{\sum_{n=0}^{\infty} \left(n+\frac{d-2}2\right)^{2s} \norm{\mathcal P_{n,d} f}_{L^2(\S^{d-1})}^2}
\\
&= \norm{\left(-\Delta^\bullet+\tfrac{(d-2)^2}{4}\right)^{s/2} f}_{L^2(\S^{d-1})}.
\label{eq:sobolev-norm-Laplacian}
\end{align}
Similarly to \eqref{eq:infty-norm-vector}, we define the Sobolev norm of a vector field $\zb f\colon\S^{d-1} \to \R^{d}$ as the Euclidean norm over its component functions $\zb f(\zb\xi) = (f_1(\zb\xi),\dots,f_{d}(\zb\xi))$, i.e.,
\begin{equation}
\norm{\zb f}_{H^{s}(\S^{d-1}\to \R^{d})}^2 = \sum_{i=1}^{d} \norm{ f_i}_{H^{s}(\S^{d-1})}^2.
\label{eq:sobolev-vector-valued}
\end{equation}
Since the negative Laplace--Beltrami operator $-\Delta^\bullet$ is self-adjoint, we can write the Sobolev norm \eqref{eq:sobolev-norm-Laplacian} as
\begin{equation*}
\norm{f}_{H^s(\S^{d-1})}^2
= \left<\left(-\Delta^\bullet+\tfrac{(d-2)^2}{4}\right)^{s} f,f \right>_{L^2(\S^{d-1})}.
\end{equation*}
We have for $s\in\N$
\begin{align*}
\norm{f}_{H^{s+1}(\S^{d-1})}^2
={}& \int_{\S^{d-1}} \left(\left(-\Delta^\bullet+\tfrac{(d-2)^2}4\right)^{s+1} f(\zb\xi)\right)  f(\zb\xi) \d\S^{d-1}(\zb\xi)
\\
={}& \int_{\S^{d-1}} \left(\left(-\Delta^\bullet+\tfrac{(d-2)^2}4\right) \left(-\Delta^\bullet+\tfrac{(d-2)^2}4\right)^s f(\zb\xi)\right) f(\zb\xi) \d\S^{d-1}(\zb\xi).
\end{align*}
Then the Green--Beltrami identity \cite[§14, Lemma 1]{Mue98}
\begin{equation}
\begin{split}
-\int_{\S^{d-1}} f(\zb\xi)\, \Delta^\bullet g(\zb\xi) \d\S^{d-1}(\zb\xi)
= \int_{\S^{d-1}} \langle\zb\nabla^\bullet f(\zb\xi),\, \zb\nabla^\bullet g(\zb\xi)\rangle \d\S^{d-1}(\zb\xi)
,\\ f\in C^{2}(\S^{d-1}),\;g\in C^{1}(\S^{d-1})
\end{split}
\label{eq:Green-Beltrami}
\end{equation}
implies that
\begin{align*}
\norm{f}_{H^{s+1}(\S^{d-1})}^2
={}& \int_{\S^{d-1}} \left< \zb\nabla^\bullet f(\zb\xi),\, \zb\nabla^\bullet \left(-\Delta^\bullet+\tfrac{(d-1)^2}4\right)^s f(\zb\xi)\right> \d\S^{d-1}(\zb\xi)
\\&+ \frac{(d-2)^2}{4} \int_{\S^{d-1}} \left(\left(-\Delta^\bullet+\tfrac{(d-2)^2}4\right)^s f(\zb\xi)\right) f(\zb\xi) \d\S^{d-1}(\zb\xi).
\end{align*}
Since the gradient $\zb\nabla^\bullet$ and the Laplacian $\Delta^\bullet$ commute by Schwarz's theorem, we obtain the recursion
\begin{equation}
\norm{f}_{H^{s+1}(\S^{d-1})}^2
= \norm{\zb\nabla^\bullet f}_{H^{s}(\S^{d-1}\to\R^{d})}^2 + \tfrac{(d-2)^2}4\norm{f}_{H^{s}(\S^{d-1})}^2.
\label{eq:sobolev-norm-gradient}
\end{equation}

\subsection{Sobolev spaces as interpolation spaces}

The norm of a bounded linear operator $\mathcal A\colon X\to Y$ between two Banach spaces $X$ and $Y$ with norms $\norm{\cdot}_X$ and $\norm{\cdot}_Y$, respectively, is defined as
\begin{equation*}
\norm{\mathcal A}_{X\to Y} = \sup_{x\in X\setminus\{0\}} \frac{\norm{\mathcal Ax}_Y}{\norm{x}_X}.
\end{equation*}

The following proposition shows that the boundedness of linear operators in Sobolev spaces $H^s(\S^2)$ can be interpolated with respect to the smoothness parameter $s$.
This result is derived from a more general interpolation theorem in \cite{Tri95}.

\begin{Proposition}
	\label{prop:sobolev-interpolation}
	Let $0\leq s_0\leq s_1$, and  
	let $\mathcal A\colon H^{s_1}(\S^{d-1})\to H^{s_1}(\S^{d-1})$ be a bounded linear operator such that its restriction to $H^{s_0}(\S^{d-1})$ is also bounded.
	For $\theta\in[0,1]$, we set $s_\theta = (1-\theta)s_0 + \theta s_1$.
	Then the restriction of $\mathcal A$ to $H^{s_\theta}(\S^{d-1})$ is bounded with
	\begin{equation*}
	\norm{\mathcal A}_{H^{s_\theta}(\S^{d-1})\to H^{s_\theta}(\S^{d-1})}
		\leq \norm{\mathcal A}_{H^{s_0}(\S^{d-1})\to H^{s_0}(\S^{d-1})}^{1-\theta} \norm{\mathcal A}_{H^{s_1}(\S^{d-1})\to H^{s_1}(\S^{d-1})}^{\theta}.
	\end{equation*}
\end{Proposition}
\begin{proof}
For $n\in\N$, let $\{Y_{n,d}^k \mid k=1,\,\dots,\,N_{n,d}\}$ be an orthonormal basis of $\mathcal P_{n,d} (L^2(\S^{d-1}))$,
and let $f\in L^2(\S^{d-1})$.
We write $f$ as the Fourier series
$$
f = \sum_{n=0}^{\infty} \sum_{k=1}^{N_{n,d}} \left<f, Y_{n,d}^k\right>_{L^2(\S^{d-1})} Y_{n,d}^k.
$$
On the index set 
$$
I=\{(n,k)\mid n\in\N,\; k=1,\dots, N_{n,d} \},
$$ 
we define for $s\ge0$ the weight function 
$$w_s(n,k) = (n+\tfrac{d-2}{2})^{2s}.$$ 
Then the Sobolev space $H^s(\S^{d-1})$ is isometrically isomorphic to the weighted $L^2$-space
\begin{equation*}
L^2 (I;w_s)
= \left\{ \hat f \colon I\to\C \middle| \norm{\smash{\hat f}}_{L^2 (I;w_s)}^2 =\sum_{(n,k)\in I} \abs{\hat f(n,k)}^2 w_s(n,k) < \infty \right\}
\end{equation*}
that consists of the Fourier coefficients 
$$\hat f(n,k)=\left<f, Y_{n,d}^k\right>_{L^2(\S^{d-1})}$$
on the set $I$ with the counting measure.
By \cite[Theorem 1.18.5]{Tri95}, the complex interpolation space between $L^2(I;w_{s_0}) \cong H^{s_0}(\S^{d-1})$ and $L^2(I;w_{s_1}) \cong H^{s_1}(\S^{d-1})$ is 
\[
\left[ L^2(I;w_{s_0}), L^2(I;w_{s_1})\right]_\theta = L^2(I;w),
\]
where 
\begin{equation*}
w(n,k) 
= \left(w_{s_0}(n,k)\right)^{1-\theta} \left(w_{s_1}(n,k)\right)^\theta 
= \left(n+\tfrac{d-2}2\right)^{2((1-\theta)s+\theta t)} 
= w_{s_{\theta}}(n,k).
\end{equation*}
Hence, $L^2(I;w) \cong H^{s_\theta}(\S^{d-1})$.
The assertion is a property of the interpolation space.
\end{proof}

\subsection{Multiplication and composition operators}
The following two theorems show that multiplication and composition with a smooth function are continuous operators in spherical Sobolev spaces $H^s(\S^{d-1})$. 

\begin{Theorem}
	\label{thm:sobolev-product}
	Let $s\in\N$.
	The multiplication operator 
	$$ H^s(\S^{d-1})\times C^s(\S^{d-1}) \to H^s(\S^{d-1}),\quad(f,v)\mapsto fv$$ 
	is continuous.
	In particular, for all $f\in H^s(\S^{d-1})$ and $v\in C^{s}(\S^{d-1})$, we have
	\begin{equation}
    \label{eq:sobolev-product}
	\norm{fv}_{H^{s}(\S^{d-1})}
	\leq c_d^s \norm{f}_{H^{s}(\S^{d-1})} \norm{v}_{C^{s}(\S^{d-1})},
	\end{equation}
	where
	$$c_d = \sqrt{2d+2}.$$
\end{Theorem}

\begin{proof}
	We use induction over $s\in\N$.
	For $s=0$, we have
	$$
	\norm{f v}_{L^2(\S^{d-1})}^2
	= \int_{\S^{d-1}} \abs{f(\zb\xi)\, v(\zb\xi)}^2 \d\zb\xi
	\leq \norm{f}_{L^2(\S^d)}^2 \norm{v}_{C(\S^{d-1})}^2.
	$$
	Let the claimed equation \eqref{eq:sobolev-product} hold for $s\in\N$, 
	and let $f\in H^{s+1}(\S^{d-1})$ and $v\in C^{s+1}(\S^{d-1})$.
	Then the decomposition \eqref{eq:sobolev-norm-gradient} of the Sobolev norm yields
	\begin{align*}
	\norm{fv}_{H^{s+1}(\S^{d-1})}^2
	&=\norm{\zb\nabla^\bullet (fv)}_{H^{s}(\S^{d-1}\to\R^d)}^2 + \tfrac{(d-2)^2}4 \norm{fv}_{H^{s}(\S^{d-1})}^2
	\\&=\norm{f\zb\nabla^\bullet v + v\zb\nabla^\bullet f}_{H^{s}(\S^{d-1}\to\R^d)}^2 + \tfrac{(d-2)^2}4 \norm{fv}_{H^{s}(\S^{d-1})}^2.
	\intertext{By the triangle inequality and since $(a+b)^2\leq2(a^2+b^2)$ for all $a,b\in\R$, we obtain}
	\norm{fv}_{H^{s+1}(\S^{d-1})}^2
	&\leq2\norm{f\zb\nabla^\bullet v}_{H^{s}(\S^{d-1}\to\R^d)}^2 + 2\norm{v\zb\nabla^\bullet f}_{H^{s}(\S^{d-1})\to\R^d}^2
	+ \tfrac{(d-2)^2}4 \norm{fv}_{H^{s}(\S^{d-1})}^2.
	\end{align*}
	By the induction hypothesis, we have
	\begin{align*}
	c_d^{-2s} \norm{f v}_{H^{d+1}(\S^{d-1})}^2
	\leq{}& 2\norm{f}_{H^{s}(\S^{d-1})}^2 \norm{\zb\nabla^\bullet v}_{C^{s}(\S^{d-1}\to\R^d)} ^2
	+ 2\norm{\zb\nabla^\bullet f}_{H^{s}(\S^{d-1}\to\R^d)}^2 \norm{v}_{C^{s}(\S^{d-1})}^2
	\\&
	+ \tfrac{(d-2)^2}4 \norm{f}_{H^{s}(\S^{d-1})}^2 \norm{v}_{C^{s}(\S^{d-1})}^2.
	\\
	={}& 2\norm{f}_{H^{s}(\S^{d-1})}^2 \norm{\zb\nabla^\bullet v}_{C^{s}(\S^{d-1}\to\R^d)} ^2
	+ \norm{\zb\nabla^\bullet f}_{H^{s}(\S^{d-1}\to\R^d)}^2 \norm{v}_{C^{s}(\S^{d-1})}^2
    \\&
	+ \norm{f}_{H^{s+1}(\S^{d-1})}^2 \norm{v}_{C^{s}(\S^{d-1})}^2 ,
	\intertext{where we made use of the decomposition  \eqref{eq:sobolev-norm-gradient} of the Sobolev norm.
    Furthermore, we apply \eqref{eq:infty-norm-gradient} and obtain}
    c_d^{-2s} \norm{f v}_{H^{d+1}(\S^{d-1})}^2
    \leq{}&  2d\norm{f}_{H^{s}(\S^{d-1})}^2 \norm{v}_{C^{s+1}(\S^{d-1})}^2
	+ \norm{f}_{H^{s+1}(\S^{d-1})}^2 \norm{v}_{C^{s}(\S^{d-1})}^2
    \\&
	+ \norm{f}_{H^{s+1}(\S^{d-1})}^2 \norm{v}_{C^{s}(\S^{d-1})}^2 .
	\end{align*}
	Because the involved norms are non-decreasing with respect to $s$, we see that
	\begin{equation*}
	\norm{f v}_{H^{s+1}(\S^{d-1})} 
	\leq c_d^s\, \sqrt{2d+2}\, \norm{f}_{H^{s+1}(\S^{d-1})} \norm{v}_{C^{s+1}(\S^{d-1})},
	\end{equation*}
    which shows \eqref{eq:sobolev-product}.
\end{proof}

\begin{Theorem}
	\label{thm:sobolev-composition}
	Let $s\in\N$, and let $\zb v\colon \S^{d-1}\to\S^{d-1}$ be bijective with $\zb v\in C^s(\S^{d-1}\to\S^{d-1})$ and $\zb v^{-1}\in C^1(\S^{d-1}\to\S^{d-1})$.
	Then there exists a constant $b_{d,s}(\zb v)$ such that for all $f\in H^s(\S^{d-1})$, we have
	\begin{equation*}
	\norm{f\circ\zb v}_{H^{s}(\S^{d-1})}
	\leq b_{d,s}(\zb v) \norm{f}_{H^{s}(\S^{d-1})}.
	\end{equation*}
\end{Theorem}
\begin{proof}
We have for $s=0$
\begin{align*}
\norm{f\circ\zb v}_{L^2(\S^{d-1})}^2
= \int_{\S^{d-1}} \abs{f(\zb v(\zb\xi))}^2 \d\S^{d-1}(\zb\xi).
\end{align*}
The substitution $\zb\eta=\zb v(\zb\xi)$ yields with the substitution rule \eqref{eq:substitution-rule}
\begin{align*}
\norm{f\circ\zb v}_{L^2(\S^{d-1})}^2
= \int_{\S^{d-1}} \abs{f(\zb\eta)}^2  \left[(\zb v^{-1})^*(\mathrm d\S^{d-1})\right] (\zb\eta).
\end{align*}
Since $\zb v^{-1}\in C^1(\S^{d-1}\to\S^{d-1})$, there exists a continuous function $\nu\colon\S^{d-1}\to\R$ such that the pullback satisfies $(\zb v^{-1})^*(\mathrm d\S^{d-1}) = \nu \d\S^{d-1}$. 
Hence, we have
\begin{equation*}
\norm{f\circ\zb v}_{L^2(\S^{d-1})}^2
\leq \norm{f}_{L^2(\S^{d-1})}^2\, \norm{\nu}_{C(\S^{d-1})},
\end{equation*}
which shows the claim for $s=0$.

We use induction on $s\in\N$.
By the decomposition \eqref{eq:sobolev-norm-gradient} of the Sobolev norm, we have
\begin{equation}
\norm{f\circ\zb v}_{H^{s+1}(\S^{d-1})}^2
= \norm{\zb\nabla^\bullet(f\circ\zb v)}_{H^{s}(\S^{d-1}\to\R^{d})}^2 +\tfrac{(d-2)^2}{4} \norm{f\circ\zb v}_{H^{s}(\S^{d-1})}^2.
\label{eq:sobolev-composition-sum}
\end{equation}
By the induction hypothesis, the second summand of \eqref{eq:sobolev-composition-sum} is bounded by
\begin{equation}
\norm{f\circ\zb v}_{H^{s}(\S^d)}
\leq b_{d,s}(\zb v) \norm{f}_{H^{s}(\S^d)}.
\label{eq:sobolev-composition-proof1}
\end{equation}
Furthermore, by \eqref{eq:sobolev-vector-valued} and the chain rule, we have for the first summand of \eqref{eq:sobolev-composition-sum}
\begin{align*}
\norm{\zb\nabla^\bullet(f\circ\zb v)}_{H^{s}(\S^{d-1}\to\R^{d})}^2
&= \norm{\zb\nabla(f\circ v)^\bullet}_{H^{s}(\S^{d-1}\to\R^{d})}^2
\\
&= \sum_{i=1}^{d} \norm{\partial_i(f\circ\zb v)^\bullet}_{H^{s}(\S^{d-1})}^2
\\&= \sum_{i=1}^{d} \norm{\sum_{j=1}^{d} ((\partial_j f^\bullet)\circ\zb v^\bullet)\, \partial_i v_j^\bullet}_{H^{s}(\S^{d-1})}^2.
\intertext{Applying the triangle inequality for the sum over $j$ and Jensen's inequality $(\sum_{j=1}^{d} x_j)^2 \leq d\,\sum_{j=1}^{d} x_j^2$, we obtain}
\norm{\zb\nabla^\bullet(f\circ\zb v)}_{H^{s}(\S^{d-1}\to\R^{d})}^2
&\leq \sum_{i=1}^{d} d \sum_{j=1}^{d} \norm{((\partial_j f^\bullet)\circ\zb v^\bullet)\, \partial_i v^\bullet_j}_{H^{s}(\S^{d-1})}^2.
\intertext{By Theorem \ref{thm:sobolev-product}, we have}
\norm{\zb\nabla^\bullet(f\circ\zb v)}_{H^{s}(\S^{d-1}\to\R^{d})}^2
&\leq  
d\, c_d^{2s} \sum_{j=1}^{d} \norm{(\partial_j f^\bullet)\circ\zb v^\bullet}_{H^{s}(\S^{d-1})}^2 
\,\sum_{i=1}^{d}\norm{\partial_i v^\bullet_j}_{C^{s}(\S^{d-1})}^2
\\
&\leq d^2\, c_d^{2s} \sum_{j=1}^{d} \norm{(\partial_j f^\bullet)\circ\zb v^\bullet}_{H^{s}(\S^{d-1})}^2 
\, \norm{v_j}_{C^{s+1}(\S^{d-1})}^2,
\end{align*}
where the last line follows from \eqref{eq:infty-norm-gradient}.
By the induction hypothesis, we see that
\begin{equation*}
\norm{\zb\nabla^\bullet(f\circ\zb v)}_{H^{s}(\S^{d-1}\to\R^{d})}^2
\leq d^2\, c_d^{2s}\, b_{d,s}(\zb v)^2\, \sum_{j=1}^{d} \norm{\partial_j f^\bullet}_{H^{s}(\S^{d-1})}^2 \,\norm{v_j}_{C^{s+1}(\S^{d-1})}^2.
\end{equation*}
By \eqref{eq:sobolev-vector-valued} and the fact that $\norm{v_j}_{C^{s+1}(\S^{d-1})}^2 \le \norm{\zb v}_{C^{s+1}(\S^{d-1}\to\R^{d})}^2$ for all $j=1,\dots,d$, we obtain
\begin{equation*}
\norm{\zb\nabla^\bullet(f\circ\zb v)}_{H^{s}(\S^{d-1}\to\R^{d})}^2
\le d^2\, c_d^{2s}\, b_{d,s}(\zb v)^2 \norm{\zb\nabla^\bullet f}_{H^{s}(\S^{d-1}\to\R^{d})}^2 \norm{\zb v}_{C^{s+1}(\S^{d-1}\to\R^{d})}^2.
\end{equation*}
Inserting the last equation and \eqref{eq:sobolev-composition-proof1} into \eqref{eq:sobolev-composition-sum}, we obtain
\begin{align*}
&\norm{f\circ\zb v}_{H^{s+1}(\S^{d-1})}^2
\\={}& \norm{\zb\nabla^\bullet(f\circ\zb v)}_{H^{s}(\S^{d-1}\to\R^{d})}^2 +\tfrac{(d-2)^2}{4} \norm{f\circ\zb v}_{H^{s}(\S^{d-1})}^2
\\
\leq{}&  b_{d,s}(\zb v)^2 \left(d^2\, c_d^{2s}\, \norm{\zb\nabla^\bullet f}_{H^{s}(\S^{d-1}\to\R^{d})}^2 \norm{\zb v}_{C^{s+1}(\S^{d-1}\to\R^{d})}^2
+ \tfrac{(d-2)^2}{4} \norm{f}_{H^s(\S^d)}^2\right)
\\
={}&  b_{d,s}(\zb v)^2 \left(\left(d^2\, c_d^{2s}  \norm{\zb v}_{C^{s+1}(\S^{d-1}\to\R^{d})}^2-1\right)\norm{\zb\nabla^\bullet f}_{H^{s}(\S^{d-1}\to\R^{d})}^2
+ \norm{f}_{H^{s+1}(\S^{d-1})}^2\right)
\\
\leq{}&  b_{d,s}(\zb v)^2 d^2\, c_d^{2s}  \norm{\zb v}_{C^{s+1}(\S^{d-1}\to\R^{d})}^2 \norm{f}_{H^{s+1}(\S^{d-1})}^2,
\end{align*}
where we have made use of \eqref{eq:sobolev-norm-gradient}.
\end{proof}

\begin{Remark}
	The last theorem resembles a similar result found in \cite[Theorem~1.2]{InKaTo13}:
	Let $M$ be a smooth, closed and oriented $d$-dimensional manifold and, for $s>\frac{d}{2}+1$, let $\varphi\in H^s(M\to M)$ be an orientation-preserving $C^1$-diffeomorphism.
	Then the composition map
	\[
	H^s(M) \to H^s(M),\;
	f \mapsto f\circ\varphi
	\]
	is continuous.
	However, \prettyref{thm:sobolev-composition} is not a special case of this result because \prettyref{thm:sobolev-composition} requires only $s\ge0$ but imposes stronger assumptions on $\varphi$.
\end{Remark}

\section{Range}

\label{sec:range}

In this section, we show that for $s\ge0$ the spherical transform
$$\mathcal U_z\colon H^s(\S^{d-1})\to H^{s+(d-2)/2}(\S^{d-1})$$
is continuous.
To this end, we seperately investigate the parts of the decomposition obtained in \prettyref{thm:spherical-transform},
$$
\mathcal U_z = \mathcal N_z \mathcal F \mathcal M_z.
$$

\subsection{Continuity of $\mathcal F$}

The following property of the Funk--Radon transform in Sobolev spaces was shown by Strichartz \cite[Lemma 4.3]{Str81} using an asymptotic analysis of its eigenvalues.

\begin{Proposition}
	\label{prop:funk-radon}
	Denote by $H^s_{\textrm{even}}(\S^{d-1})$ the restriction of the Sobolev space $H^s(\S^{d-1})$ to even functions $f(\zb\xi)=f(-\zb\xi)$.
	The Funk--Radon transform
	\begin{equation*}
	\mathcal F\colon H^s_{\textrm{even}}(\S^{d-1})\to H_{\textrm{even}}^{s+\frac{d-2}{2}}(\S^{d-1})
	\end{equation*}
	is continuous and bijective.
\end{Proposition}

\subsection{Continuity of $\mathcal M_z$ and $\mathcal N_z$}

\begin{Theorem}
	\label{thm:mn}
	Let $z\in[0,1)$ and $s \in\R$ with $s\ge0$.
	The operators 
	$$\mathcal M_z\colon H^s(\S^{d-1})\to H^s(\S^{d-1})$$
	and 
	$$\mathcal N_z\colon H^s(\S^{d-1})\to H^s(\S^{d-1}),$$ 
	as defined in \eqref{eq:M} and \eqref{eq:N}, are continuous and open.
\end{Theorem}
\begin{proof}
	We first perform the proof for $\mathcal M_z$.
	Initially, we consider only the situation $s\in\N$.
	Let $f\in H^s(\S^{d-1})$ and $z\in(-1,1)$.
	We write
	\begin{equation*}
	\mathcal M_z f (\zb{\xi}) =
	w_z(\zb\xi)\,
	[f \circ \zb h_z] ( \zb \xi),
	\qquad\zb\xi\in\S^{d-1},
	\end{equation*}
	where
	\begin{equation*}
	w_z\colon \S^{d-1}\to\R,\quad
	w_z(\zb\xi)
	= \left(\frac{\sqrt{1-z^2}}{1+z\xi_{d}}\right)^{d-2}
	\end{equation*}
	and $\zb h_z$ is given in \eqref{eq:h_z}.
	We see that the extension
	\begin{equation*}
	w_z^\bullet(\zb x)
	= w_z \left(\frac{\zb x}{\norm{\zb x}}\right)
	= \left(\sqrt{1-z^2}\,\frac{\norm{\zb x}}{\norm{\zb x}+zx_{d}}\right)^{d-2}
    ,\qquad\zb x\in\R^d\setminus\{\zb0\},
	\end{equation*}
	is smooth except in the origin,
	i.e., $w_z^\bullet \in C^\infty(\R^{d}\setminus\{\zb0\})$.
	Hence, $w_z\in C^\infty(\S^{d-1})$.
	Then Theorem \ref{thm:sobolev-product} implies that
	\begin{equation}
	\norm{w_z\,(f\circ \zb h_z)}_{H^s(\S^{d-1})}
	\le c_{d}^s\, \norm{w_z}_{C^s(\S^{d-1})} \norm{f\circ \zb h_z}_{H^s(\S^{d-1})}.
    \label{eq:norm-w-fh}
	\end{equation}
	Moreover, the extension of $\zb h_z$,
	\begin{equation*}
	\zb h_z^\bullet \colon \R^{d}\setminus\{\zb0\}\to\R^{d},\quad
	\zb h_z^\bullet (\zb x)
	= \sum_{i=1}^{d-1} \frac{\sqrt{1-z^2}}{\norm{\zb x}+zx_{d}} x_i \zb \epsilon^i + \frac{z\norm{\zb x}+x_{d}}{\norm{\zb x}+zx_{d}} \zb \epsilon^{d},
	\end{equation*}
	is also smooth, so $\zb h_z\in C^\infty(\S^{d-1}\to\S^{d-1})$.
	This implies that also the inverse
	$\zb h_z^{-1} = \zb h_{-z}$, see \eqref{eq:h^-1}, is smooth.
	So $\zb h_z$ is a diffeomorphism and Theorem \ref{thm:sobolev-composition} together with \eqref{eq:norm-w-fh} implies that
	\begin{align*}
	\norm{\mathcal M_zf}_{H^s(\S^{d-1})}
    &\le c_{d}^s\, \norm{w_z}_{C^s(\S^{d-1})} \norm{f\circ \zb h_z}_{H^s(\S^{d-1})}
    \\
    &\le c_{d}^s\, \norm{w_z}_{C^s(\S^{d-1})} b_{d,s}(\zb h_z) \norm{f}_{H^s(\S^{d-1})}.
	\end{align*}
	Thus, the operator $\mathcal M_z\colon H^s(\S^{d-1})\to H^s(\S^{d-1})$ is continuous.
	
	Now let $s\in\R$ with $s\ge0$.
	The above proof shows that both the restrictions of $\mathcal M_z$ to $H^{\lfloor s\rfloor}$ and to $H^{\lfloor s\rfloor +1}$ are continuous, where $\lfloor s\rfloor$ denotes the largest integer that is smaller than or equal to $s$.
	The continuity of $\mathcal M_z$ on $H^s(\S^{d-1})$ follows by the interpolation result \prettyref{prop:sobolev-interpolation}.
	
	In order to prove the openness of $\mathcal M_z$, we show that the inverse $\mathcal M_z^{-1}$ restricted to $H^s(\S^{d-1})$ is continuous.
    However, we have already done this because $\mathcal M_z^{-1} = \mathcal M_{-z}$ by \eqref{eq:M-1}.
	
	The same argumentation as above also works for the operator $\mathcal N_z$ as follows.
    Let $z\in[0,1)$ and $s\in\N$.
	We write
	$$\mathcal N_zf (\zb{\xi}) =
	v_z(\zb\xi)\,
	[f \circ \zb g_z] ( \zb \xi)
    ,\qquad\zb\xi\in\S^{d-1},
	$$
	where 
    $$v_z(\zb\xi) =  (1-z^2\xi_{d}^2)^{-\frac{d-2}{2}}
    ,\qquad\zb\xi\in\S^{d-1}.$$
	We see that the extension
	$$
	v_z^\bullet(\zb x) =  \left(\frac{\norm{\zb x}^2}{\norm{\zb x}^2-z^2x_{d}^2}\right)^{\frac{d-2}{2}},
	\qquad\zb x\in\R^{d}\setminus\{\zb0\},
	$$ 
	is smooth and hence $v_z\in C^s(\S^d)$.
	\prettyref{thm:sobolev-product} yields
	\begin{equation*}
	\norm{v_z\, (f\circ \zb g_z)}_{H^{s}(\S^{d-1})} 
	\leq c_d^s \norm{v_z}_{C^s(\S^{d-1})}  \norm{f\circ \zb g_z}_{H^{s}(\S^{d-1})}.
	\end{equation*}
	Since the extensions of both
	\begin{equation*}
	\zb g_z^\bullet(\zb x)
	=\zb g_z \left(\frac{\zb x}{\norm{\zb x}}\right)
	= \frac{1}{\sqrt{\norm{\zb x}^2-z^2{x_{d}^2}}} 
	\left(\sum_{i=1}^{d-1} {x_i} \zb \epsilon^i + \sqrt{1-z^2}\, {x_{d}} \zb \epsilon^{d}\right)
    , \qquad \zb x\in\R^d\setminus\{\zb0\},
	\end{equation*}
	and its inverse \eqref{eq:g^-1}
	\[
	[(\zb g_z^{-1})^\bullet](\zb x)
	=\frac{1}{\sqrt{\norm{\zb x}^2-z^2+z^2x_{d}^2}} \left(\sqrt{1-z^2}\,\sum_{i=1}^{d-1} \omega_i \zb\epsilon^i
	+x_{d} \zb\epsilon^{d}\right)
    , \qquad \zb x\in\R^d\setminus\{\zb0\},
	\]
	are smooth functions on $\R^{d}\setminus\{\zb 0\}$, we see that $\zb g_z$ is a smooth diffeomorphism in $C^s(\S^{d-1})$.
	By \prettyref{thm:sobolev-composition}, there exists a constant $b_{d,s} (\zb g_z)$ independent of $f$ such that 
	$$\norm{f\circ \zb g_z}_{H^{s}(\S^{d-1})} \leq b_{d,s} (\zb g_z) \norm{f}_{H^{s}(\S^{d-1})}.$$
	Hence, $\mathcal N_z$ is a bounded operator on $H^s(\S^{d-1})$.
	An analogue computation shows that the inverse operator \eqref{eq:N-1}
	\begin{equation*}
	\mathcal N_z^{-1} f (\zb\eta)
	= \frac{f(\zb g_z^{-1}(\zb\eta))}{v_z(\zb g_z^{-1}(\zb\eta))}
	= \left(\frac{1-z^2}{1-z^2+z^2\eta_{d}^2}\right)^{\frac{d-2}{2}}\,(f\circ \zb g_z^{-1})(\zb\eta)
    ,\qquad\zb\eta\in\S^{d-1}
	\end{equation*}
	is also bounded on $H^s(\S^{d-1})$.
    The assertion for general $s$ follows by the same interpolation argument as for $\mathcal M_z$.
\end{proof}

\subsection{Continuity of $\mathcal U_z$}

\begin{Theorem}
	\label{thm:Uz-Sobolev}
	Let $z\in(0,1)$ and $s\in\R$ with $s\ge0$.
	We set $H^s_z(\S^{d-1})$ as the subspace of all functions $f\in H^s(\S^{d-1})$ that satisfy
	\begin{equation}
    f(\zb{\omega})
	=\left(\frac{1-z^2}{1-2z\omega_{d}+z^2}\right)^{d-2} f \circ\zb r_z(\zb\omega),
	\qquad\zb\omega\in\S^{d-1},
    \label{eq:symmetric-z}
    \end{equation}
    almost everywhere,
	where the point reflection $\zb r_z$ about the point $z\:\!\zb\epsilon^d$ is given in \eqref{eq:point-reflection}.
	Then the spherical transform
	\begin{equation*}
	\mathcal U_z\colon H^s_z(\S^{d-1})\to H^{s+\frac{d-2}{2}}_{\textrm{even}}(\S^{d-1})
	\end{equation*}
	is continuous and bijective and its inverse operator is also continuous.
\end{Theorem}
\begin{proof} In \prettyref{thm:spherical-transform}, we obtained the decomposition
	$$
	\mathcal U_z = \mathcal N_z \mathcal F \mathcal M_z.
	$$
	We are going to look at the parts of this decomposition separately.
	By Theorem \ref{thm:mn}, we obtain that
	\begin{equation*}
	\mathcal M_z \colon H^s(\S^{d-1})\to H^s(\S^{d-1})
	\end{equation*}
	is continuous and bijective.
	The same holds for the restriction
	\begin{equation*}
	\mathcal M_z \colon H^s_{z} (\S^{d-1})\to H^s_{\textrm{even}}(\S^{d-1}),
	\end{equation*}
	which follows from the characterization of the nullspace in Theorem \ref{thm:nullspace}.
	By Proposition \ref{prop:funk-radon}, the Funk--Radon transform
	\begin{equation*}
	\mathcal F \colon H^s_{\textrm{even}}(\S^{d-1})\to H^{s+\frac{d-2}{2}}_{\textrm{even}}(\S^{d-1})
	\end{equation*}
	is continuous and bijective.
	Finally, Theorem \ref{thm:mn} and the observation that any function $f\colon\S^{d-1}\to\C$ is even if and only if $\mathcal N_zf$ is even show that
	\begin{equation*}
	\mathcal N_z\colon H^{s+\frac{d-2}{2}}_{\textrm{even}}(\S^{d-1})\to H^{s+\frac{d-2}{2}}_{\textrm{even}}(\S^{d-1})
	\end{equation*}
	is continuous and bijective.
	The continuity of the inverse operator of $\mathcal U_z$ follows from the open mapping theorem.
\end{proof}

Theorem \ref{thm:Uz-Sobolev} is a generalization of \prettyref{prop:funk-radon} for the Funk--Radon transform $\mathcal F$; the main difference is that the space $H^s_{\textrm{even}}(\S^{d-1})$ is replaced by $H^s_z(\S^{d-1})$, which contains functions that satisfy the symmetry condition \eqref{eq:symmetric-z} with respect to the point reflection in $z\,\zb\epsilon^d$.
Furthermore, the spherical transform $\mathcal U_z$ is smoothing of degree $\frac{d-2}{2}$, which comes from the fact that $\mathcal U_z$ takes the integrals along $(d-2)$-dimensional submanifolds.

\section{Geometric interpretation}
\label{sec:geometric-interpratation}

We give geometric interpretations of the mappings $\zb g_z$ and $\zb h_z\colon \S^{d-1}\to\S^{d-1}$ that were defined in \prettyref{sec:hg}.
The mapping $\zb g_z$ consists of a scaling with the factor $\sqrt{1-z^2}$ along the $\xi_{d}$ axis, which maps the sphere to an ellipsoid, which is symmetric with respect to rotations about the $\xi_{d}$ axis.
Then a central projection maps this ellipsoid onto the sphere again.

In order give a description of the mapping $\zb h_z$, we define the stereographic projection
$$
\zb\pi\colon \S^{d-1}\setminus\{\zb\epsilon^{d}\}\to\R^{d-1},\quad
\zb\xi \mapsto \sum_{i=1}^{d-1} \frac{\xi_i}{1-\xi_{d}} \zb\epsilon^i
$$
and its inverse
$$
\zb\pi^{-1}\colon \R^{d-1}\to \S^{d-1}\setminus\{\zb\epsilon^{d}\},\quad
\zb x \mapsto \frac{2\zb x + (\norm{\zb x}^2-1)\, \zb\epsilon^{d}}{1+\norm{\zb x}^2}.
$$
The following corollary states that, via stereographic projection, the map $\zb h_z$ on the sphere $\S^{d-1}$ corresponds to a uniform scaling in the equatorial hyperplane $\R^{d-1}$ with the scaling factor $\sqrt{\frac{1+z}{1-z}}$.

\begin{Corollary}
	Let $\zb\xi\in\S^{d-1}$ and $z\in(0,1)$. Then we have
	\begin{equation*}
	\zb h_z(\zb\xi) =
	\zb\pi^{-1}\left(\sqrt{\frac{1+z}{1-z}} \,\zb\pi(\zb\xi) \right).
	\end{equation*}
\end{Corollary}
\begin{proof}
    We are going to show that
    \begin{equation*}
    \zb\pi(\zb h_z(\zb\xi))
    = \sqrt{\frac{1+z}{1-z}}\,\zb\pi(\zb\xi)
    \end{equation*}
    holds.
	We have on the one hand
	$$
	\sqrt{\frac{1+z}{1-z}}\,\zb\pi(\zb\xi)
	= \sqrt{\frac{1+z}{1-z}}\, \sum_{i=1}^{d-1} \frac{\xi_i}{1-\xi_{d}} \zb\epsilon^i
	$$
    and the other hand
    \begin{align*}
    \zb\pi(\zb h_z(\zb\xi))
    &= \sum_{i=1}^{d-1} \frac{\frac{\sqrt{1-z^2}\,\xi_i}{1+z\xi_d}}{1-\frac{z+\xi_d}{1+z\xi_d}}\,\zb\epsilon^i
    = \sum_{i=1}^{d-1} \frac{\sqrt{1-z^2}\,\xi_i}{1+z\xi_d-(z+\xi_d)}\,\zb\epsilon^i
    \\&
    = \sum_{i=1}^{d-1} \frac{\sqrt{1-z^2}\,\xi_i}{(1-z)\,(1-\xi_d)}\,\zb\epsilon^i.
    \end{align*}
    The assertion follows by canceling $\sqrt{1-z}$.
\end{proof}

\end{document}